\documentclass[a4paper,11pt,twoside,tbtags]{article} 
\input{TT-arxive.sty}

\DN\Jli{J _{\ell ,\ii }}
\DN\Jonei{J _{1,\ii }}
\DN\Jtwoi{J _{2,\ii }}
\DN\rfrak{r}

\DN\iii{\mathbf{i}}
\DN\jjj{\mathbf{j}}
\DN\iiim{(\mathbf{i}_1,\ldots,\mathbf{i}_m)}
\DN\iiM{\ii }

\DN\Fl{\mathcal{F}_{\ell } }
\DN\ulab{\mathfrak{u}}
\DN\ulabl{\mathfrak{u}_{\ell }}
\DN\ulabi{\mathfrak{u}_{\ell ,\mathsf{i}}}

\DN\FFF{\mathbb{F}} \DN\FFFl{\FFF (\ell )}\DN\FFFz{\FFF (0)}
\DN\Il{I (\ell )}

\DN\mmFliiis{\mmm _{\mm _{\ell , \iii _{\sigma }} }}

\DN\wIIIl{\wIII (\ell )} 
\DN\wIII{\widetilde{\III }}
\DN\mmFl{\mmm _{\FFFl }}

\DN\I{I}
\DN\lambdal{\lambda _{\IIIl }}
\DN\III{\mathbb{I}}
\DN\IIIxl{\III (\ell )_{\mathbf{x}}} 
\DN\IIIl{\III (\ell )} 
\DN\IIl{\mathsf{I}(\ell )}

\DN\AAAA{\mathcal{A}}
\DN\Ai{\mathcal{A}_i } \DN\Aj{\mathcal{A}_j }
\DN\Aii{\mathcal{A}_{i'} } \DN\Ajj{\mathcal{A}_{j'} }
\DN\Aiii{\mathcal{A}_{i^-} } \DN\Ajjj{\mathcal{A}_{j^-} }
 \DN\Aiin{\mathcal{A}_{\ii _n}} 
 \DN\Ali{\AAAA _{\ell , \ii } } 
 \DN\Aliii{\AAAA _{\ell , \ii ^-} } 
 \DN\Azi{\AAAA _{0 , \ii } } 
\DN\iim{(\ii _1,\ldots,\ii _m)}

\DN\B{\mathcal{B} }
\DN\Bi{\B _{\ii }} 
\DN\Bl{\B _{\ell }}
\DN\Bli{\B _{\ell , \ii }}
\DN\Blin{\B _{\ell , \ii _n}}
\DN\Blj{\B _{\ell , \jj }}
\DN\Bj{\B _{\jj }}

\DN\Gl{\mathcal{G}_{\ell }} \DN\Gll{\mathcal{G}_{\ell +1}}

\DN\fliirs{f_{\ell ,\ii _{r,s}}}
\DN\fliir{f_{\ell ,\ii _r}}
\DN\flii{f_{\ell ,\ii }}\DN\fljj{f_{\ell , \jj }}
\DN\fliip{f_{\ell ,\iip }}
\DN\fliiq{f_{\ell ,\iiq }}
\DN\fliipsigma{f_{\ell ,\ii _{\sigma (p )} }}
\DN\fljjq{f_{\ell , \jjq }}
\DN\fljjp{f_{\ell , \jjp }}

\DN\fljjr{f_{\ell , \jj _r}}
\DN\fljjrs{f_{\ell , \jj _{r,s}}}
\DN\fzii{f_{0 ,\ii }}\DN\fzjj{f_{0 , \jj }}

\DN\fii{f_{\ii }}\DN\fjj{f_{\jj }}
\DN\fiil{f_{\ii '}}
\DN\fliil{f_{\ell , \ii '}}
\DN\ii{i}\DN\jj{j}

\DN\iip{\ii _p}
\DN\jjp{\jj _p}
\DN\jjq{\jj _q}
\DN\iiq{\ii _q}

\DN\rhoF{\rho _{\mathbb{F}}}
\DN\rhoFl{\rho _{\FFFl }}

\DN\rhoZ{\rho _{\Gz }}
\DN\rhol{\rho _{\Gl }}
\DN\nuF{\nu_{\mathbb{F}}}
\DN\nuFl{\nu_{\FFFl }}
\DN\KF{\K _{\mathbb{F}}}
\DN\KFl{\K _{\FFFl }} 
\DN\f{\mathbf{f}}
\DN\Det{\mathrm{Det}}

\DN\aaa{\mathsf{a}}
\DN\ijn{{_{i,j=1}^{n}}}
\DN\mub{\mu _{\mathrm{Be},\alpha }}
\DN\rgn{\rho _{\mathrm{Gin}}}
\DN\kg{\mathsf{K}_{\mathrm{Gin}}}

\DN\mug{\mu _{\mathrm{Gin}}}
\DN\mugx{\mu _{\mathrm{Gin},x}}
\DN\muSinb{\muone _{\mathrm{Sin}, \beta }}
\DN\muAi{\mu _{\mathrm{Ai}}}
\DN\lab{\mathfrak{l}}


\DN\rairybeta{\rho _{\mathrm{Ai},\beta }}
\DN\rairytwo{\rho _{\mathrm{Ai}}} 
\DN\rNone{\rho ^{\n ,1}}
\DN\rNtwo{\rho ^{\n ,2}}
\DN\rNnk{\rho ^{\n ,n+k}}

\DN\rhobar{\bar{\rho }}
\DN\rb{\rho _{\beta }^n}
\DN\rbone{\rho _{\beta }^1}
\DN\rbx{\rho _{1 }^n}
\DN\rby{\rho _{2 }^n}
\DN\rbz{\rho _{4 }^n}

\DN\rbNone{\rho _{\beta }^{\n ,1}}
\DN\rbNtwo{\rho _{\beta }^{\n ,2}}
\DN\rbNxone{\rho _{\beta ,x}^{\n ,1}}
\DN\rbNxtwo{\rho _{\beta ,x}^{\n ,2}}
\DN\rbNxn{\rho _{\beta ,x}^{\n ,n}}
\DN\rN{\rho ^{\n ,n}}
\DN\rbN{\rho _{\beta }^{\n ,n}}
\DN\rbNN{\rho _{\beta }^{\n ,N}}
\DN\rbNx{\rho _{1}^{\n ,n}}
\DN\rbNy{\rho _{2}^{\n ,n}}
\DN\rbNz{\rho _{4}^{\n ,n}}


\DN\vxNi{v _{\xNi }} 	
\DN\vyNi{v _{\yNi }}
\DN\GN{G_{\nN }} 	
\DN\GNN{G_{\nN }\ts \N } 	
\DN\CGN{\C ^{\GN\ts \N }}
\DN\CN{\C ^{\N }} 	\DN\RN{\R ^{\N }}
\DN\xyGN{[x]_{\nN }, \, [y]_{\nN } \in \GN }
\DN\xGN{[x]_{\nN }  \in \GN }
\DN\yGN{[y]_{\nN }  \in \GN }

\DN\vv{\mathbf{v}}
\DN\axyN{\mathsf{K}_{x,y}^{\nN }}
\DN\MMN{\mathsf{M}^{\nN }}
\DN\MN{\mathsf{M}^{\nN }}
\DN\LLN{\LL ^{\nN }}
\DN\LLNx{\LLN _{x}}
\DN\LLNy{\LLN _{y}}
\DN\LLxNi{\LL _{\xNi }}
\DN\LLyNi{\LL _{\yNi }}
\DN\LLxi{\LL _{\xi }}
\DN\LLeta{\LL _{\eta }}

\DN\T{\mathbf{T}}
\DN\Tl{\T _{\ell }}
\DN\Tll{\T _{\ell '}}
\DN\TTT{\mathbb{T}}
\DN\Tail{\mathrm{Tail} }
\DN\nN{\mathsf{N}}
\DN\IN{I_{\nN }}
\DN\sS{\mathrm{S}}
\DN\sSS{\mathfrak{S}}
\DN\sss{\mathsf{s}}
\DN\K{\mathsf{K}}
\DN\Kl{{\mathsf{K}}_{\ell }}
\DN\KKK{\mathsf{K}} 
\DN\LL{\mathsf{L}}
\DN\mm{ f }
\DN\mmm{\mathsf{m}}
\DN\xN{[x]_{\nN }}
\DN\yN{[y]_{\nN }}

\DN\TailS{\Tail (\Ssf )}
\DN\Path{\mathsf{Path}}
\DN\Tree{\mathsf{Tree}}
\DN\Da{\underline{\Delta }}\DN\Db{\underline{\Delta }'}
\DN\DDa{\{\Dl \}_{\ell \in \{ 0 \}\cup \N} }

\DN\Ssf{\mathsf{S}}
\DN\rrr{\mathsf{r}}
\DN\Tf{\Tree _{\mathrm{fin}}}

\DN\rr{\mathsf{r}}
\DN\eN{\ell \in \N }
\DN\nnn{\mathsf{n}}
\DN\nNN{\nN _{\nnn }}
\DN\muN{\mu ^{\nnn }}

\DN\li{_{\ell }(i)}
\DN\kl{_{\ell }}
\DN\kll{_{\ell +1}}
\DN\kkl{_{\ell }}
\DN\Kli{\K _{\ell ,\infty}} 
\DN\Ki{\K _{\infty }}
\DN\Klbar{\K  \kl }
\DN\KDkl {\K \kl }
\DN\ml{\mm \kl } 
\DN\mlv{\mm _{\ell , v} } 
\DN\mli{\mm _{\ell }(i)} \DN\mlj{\mm _{\ell }(j)}

\DN\ff{\mathfrak{f}}

\DN\Klt{\widetilde{\K }_{\ell }}

\DN\Dz{\Delta (0)}\DN\Done{\Delta (1)}
\DN\Dl{\Delta (\ell ) } \DN\Dll{\Delta (\ell +1) }  
\DN\Di{\Delta _{\infty }  } 
\DN\n{\mathsf{n}}
\DN\nd{\n \Delta }
\DN\hatx{\hat{x}}
\DN\Sr{\Sit _r}
\DN\Sit{S}
 
\DN\mmF{\mmm _{\mm _{\ell }}}  

\DN\mmFiii{\mmm _{\mm _{\ell ,\iii }}} 
\DN\mmFiiii{\mmm _{\mm _{\ell , (\ii _1,\ldots,\ii _m)}} }
\DN\mmFiiin{\mmm _{\mm _{\ell , \iii _n}} }

\DN\mmFi{\mmm _{\mm _{\ell , \iiM }}}
\DN\mmFii{\mmm_{\mm _{\ell , \mathsf{i}}} }

\DN\muFi{\mu _{\FFFl ,\iiM }}
\DN\muFiii{\mu _{\FFFl ,\iii }}
\DN\muk{\mu }
\DN\muz{\mu _{\Dz }}

\DN\mul{\mu |_{\Gl }} 
\DN\mulsA{\muk ( \mathsf{A} | \Gl  )(\mathsf{s}) }
\DN\mulA{\muk ( \mathsf{A} | \Gl  )} 
\DN\mulsU{\muk ( U | \Gl  )(\mathsf{s}) }

\DN\muls{\muk (\cdot | \Gl  )}
\DN\muF{\mu _{\mathbb{F}}}
\DN\nut{\widetilde{\nu }}
\DN\nul{\nuFl } 
\DN\nuli{\nu _{\ell ,\infty } } 
\DN\nulit{\nut _{\ell ,\infty } } 
\DN\nult{\nut _{\ell }}
\DN\nulc{\check{\nu }_{\ell }}

\DN\SSi{\Ssf _{\infty }}
\DN\Omegat{\widetilde{\Omega}}
\DN\Omegac{\check{\Omega }} \DN\Omegacl{\check{\Omega }_{\ell }}
\DN\tauc{\check{\tau }}
\DN\FF{\mathsf{F}}

\DN\kappal{\kappa _{\ell ,\mathsf{i}}}
\DN\iotal{\iota_{\ell} }
\DN\etai{\eta _{\infty} }
\DN\iotai{\iota _{\infty} }

 \DN\dist{\mathsf{d}}

\DN\0{\rhoFl ^m (\iii )}
\DN\1{\Omega }
\DN\2{\mathrm{Conf}(\1 )}
\DN\3{\overline{\Omega}(\ell )}
\DN\4{\underline{\Omega}(\ell )}
\DN\5{\mu _{\FFFl ,\mathsf{i}}}
\DN\6{\FFFl = \{ \flii \}_{\ii \in\IIIl }}
\DN\7{ \nul \diamond \mmFl }
\DN\8{ ( \nul \diamond \mmF ) \circ \ulabl ^{-1}|_{\Gl } }
\DN\9{\prod_{n=1}^m \Blin }

\title{
Discrete approximations of determinantal point processes 
on continuous spaces: tree representations and  tail triviality}
\author{Hirofumi \textsc{Osada}$^1 $, Shota \textsc{Osada}$ ^2$ \\
{\small $ ^1$osada@math.kyushu-u.ac.jp,$ ^2$s\_osada@kyudai.jp}
} \maketitle


\abstract{We prove  tail triviality of determinantal point processes 
$ \mu $ on continuous spaces. 
Tail triviality had been proved for such processes only on 
discrete spaces, and hence we have generalized the result to continuous spaces. 
To do this, we 
construct tree representations, that is, discrete approximations of determinantal point processes enjoying a determinantal structure. 
There are many interesting examples of determinantal point processes on continuous spaces such as zero points of the hyperbolic Gaussian analytic function with Bergman kernel, and the thermodynamic limit of eigenvalues of Gaussian random matrices for Sine$_2 $, Airy$_2 $, Bessel$_2 $, and Ginibre point processes. 
Tail triviality of $ \mu $ plays a significant role in proving 
the uniqueness of solutions of infinite-dimensional stochastic differential equations (ISDEs) associated with $ \mu $. 
For particle systems in $ \R $ arising from random matrix theory, 
there are two completely different constructions of natural stochastic dynamics. 
One is given by stochastic analysis through ISDEs and Dirichlet form theory, and the other is an algebraic method based on space-time correlation functions. 
Tail triviality is used crucially to prove the equivalence of these two stochastic dynamics. 
}


\DN\proofbegin{\smallskip \noindent {\it Proof. }}
\DN\proofend{\qed \smallskip }

\section{Introduction}\label{s:1}

Let $ \Sit $ be a locally compact, complete, separable metric space with 
 metric $ \dist (\cdot, \cdot)$.  
 We assume $\Sit $ is unbounded. 
We equip $ \Sit $ with a Radon measure 
$ \mmm $ such that $ \mmm (\mathcal{O} ) > 0 $ 
for any non-empty open set $ \mathcal{O} $ in $ \Sit $. 
Let $ \Ssf $ be the configuration space over $ \Sit $ 
(see \eqref{:20a} for definition). 
$ \Ssf $ is a Polish space equipped with the vague topology.

A determinantal point process $ \mu $ on $ \Sit $ 
is a probability measure on $ (\Ssf ,\mathcal{B}(\Ssf ) )$ 
for which the $ m $-point correlation function $ \rho ^m $ 
with respect to $ \mmm $ is given by the determinant 
\begin{align}\label{:10a}&
\rho ^m (\mathbf{x}) = \det [\mathsf{K}(x_i,x_j)]_{i,j=1}^m 
.\end{align}
Here $ \map{\K }{\Sit \ts \Sit }{\mathbb{C}}$ is a measurable kernel and 
$ \mathbf{x} = (x_1,\ldots,x_m )$. 
We refer to \sref{s:2} and {\it e.g.} \cite{AGZ10,forrester,sosh} for the definition of correlation functions and related notions. 
$ \mu $ is said to be 
associated with $ (\mathsf{K}, \mmm )$ and also a 
$ (\mathsf{K}, \mmm )$-determinantal point process.

We set $ \KKK f (x)= \int_{\Sit } \K (x,y) f (y) \mmm (dy)$. 
We regard $ \K $ as an operator on $ L^2 (\Sit ,\mmm )$ and denote it by the same symbol. 
We say $ \KKK $ is of locally trace class if 
$ \KKK _{A} f (x)= \int 1_{A} (x) \K (x,y) 1_{A} (y)f(y)\mmm (dy)$ 
is a trace class operator on $ L^2(\Sit , \mmm )$  for any compact set $ A $.

Throughout this paper, we assume that $ \K $ satisfies:

\smallskip 

\noindent 
\thetag{A1} $ \K $ is bounded, Hermitian symmetric, 
of locally trace class, and $ \mathrm{Spec}(\KKK  )\subset [0,1]$. 

\smallskip 
%

From \thetag{A1} we deduce that 
 the associated determinantal point process $ \mu = \mu ^{\K ,\mmm }$ 
exists and is unique \cite{sosh,s-t.aop,RL.03}.

In the last two decades, determinantal point processes have been extensively studied. They contain many interesting examples; e.g., spanning trees and Schur measures on discrete spaces, zero points of the hyperbolic Gaussian analytic function with Bergman kernel, and thermodynamic limits of eigenvalues of Gaussian random matrices such as Sine$_2 $, Airy$_2 $, Bessel$_2 $, and Ginibre point processes on continuous spaces \cite{AGZ10,hkpv.gaf,sosh}.

Determinantal point processes on discrete spaces have a well-behaved algebraic structure; as a result, some important facts are only known for discrete determinantal point processes \cite{ghosh.14,RL.03duke,RL.03,RL.14,s-t.aop}. 
One such example is  tail triviality, which says that each event of a tail $ \sigma $-field $ \TailS $ takes value 0 or 1.  We refer to \eqref{:20c} for the definition of 
$ \TailS $. 

The purpose of this paper is to prove that the tail $ \sigma $-field $ \TailS  $ of $ \Ssf $ is trivial with respect to $ \mu $. 
If the space $ \Sit $ is discrete, then  tail triviality has been proved by 
Shirai-Takahashi \cite{s-t.jfa} for 
$ \mathrm{Spec}(\K )\subset (0,1)$, and by Russell Lyons \cite{RL.03} 
for $ \mathrm{Spec}(\K )\subset [0,1]$. 
If the space $ \Sit $ is continuous, the problem remained open \cite{RL.14}. 

To prove  tail triviality we introduce a discrete approximation for determinantal point processes, called the tree representation. This representation has a determinantal structure, and so belongs to determinantal point processes on discrete spaces. 

A $ \mmm $-partition $ \Delta =\{ \mathcal{A}_i  \}_{i\in I } $ of $ \Sit $ is a countable collection of disjoint relatively compact, measurable subsets of $ \Sit $ such that $ \cup_i \mathcal{A}_i = \Sit  $ and that 
$ \mmm ( \mathcal{A}_i ) > 0 $ for all $ i \in I $. 
For two partitions $ \Delta =\{ \mathcal{A}_i  \}_{i\in I } $ and 
$ \Gamma=\{ \mathcal{B}_{j}  \}_{j\in J } $, we write 
$ \Delta \prec \Gamma $ if for each $ j \in J $ there exists $ i \in I $ 
such that $ \mathcal{B}_{j}  \subset \mathcal{A}_i $. 
We assume: 

\medskip

\noindent 
\thetag{A2} There exists a sequence of $ \mmm $-partitions 
$ \{ \Delta (\ell ) \}_{\ell \in \N } $ satisfying 
\eqref{:10b}--\eqref{:10d}. 
\begin{align}  \label{:10b}&\quad \quad \quad 
\Dl \prec \Dll \quad \text{ for all } \eN 
,\\ 
\label{:10c}&\quad \quad \quad 
\sigma [\bigcup_{\eN }\mathcal{F}\kl ] = \mathcal{B}(\Sit ) 
,
\\\label{:10d}& \quad \quad \quad 
\# \{ \jj ; \AAAA _{\ell + 1,\jj } \subset \Ali \} = 2 
\text{ for all } \ii \in \Il \text{ and } \ell \in \N 
,\end{align}
where we set 
$ \Dl = \{ \Ali  \}_{i\in \Il } $ and 
$ \mathcal{F}\kl := \mathcal{F}_{\Dl } = \sigma [ \Ali ; i\in \Il  ]$. 
Furthermore, $ \# \{ \cdot  \}  $ denotes the cardinality of $  \{\cdot \}$. 
\ms 

Condition \eqref{:10d} is just for simplicity. 
This condition implies that the sequence $ \{ \Delta (\ell ) \}_{\ell \in \N } $ 
has a binary tree-like structure. 
We remark that \thetag{A2} is a mild assumption and, indeed, 
satisfied if $ \Sit $ is an open set in $ \Rd $ 
and $ \mmm $ has positive density with respect to the Lebesgue measure. 
We now state one of our main theorems: 
\begin{theorem}	\label{l:11} 
Assume \thetag{A1} and \thetag{A2}. 
Let $ \mu $ be the $ (\K , \mmm )$-determinantal point process. 
Then $ \mu $ has a trivial tail. That is, 
$ \mu (A) \in \{ 0,1 \} $  for all $ A \in \TailS $. 
\end{theorem}

Many interesting determinantal point processes arise from random matrices 
such as  Sine$_2 $, Airy$_2 $, and Bessel$_2 $ point processes 
in $ \mathbb{R}$ and 
the Ginibre point process in $ \mathbb{R}^2$. 
Applying \tref{l:11} to these examples we obtain that all have trivial tails. 
We shall present these examples in \sref{s:7}. 

%

We now explain the idea of the proof. 
We have two candidates for the discrete approximations of $ \mu $. 
One is the approximation of the kernel $ \K $. Let $ \K _{\ell } (x,y)$ be the 
discrete kernel on $ \Il $ such that 
\begin{align}\notag & 
\K _{\ell } (x,y) = 
\frac{1}{\mmm (\mathcal{A}_{\ell }(x))\mmm (\mathcal{A}_{\ell }(y))} 
\int_{\mathcal{A}_{\ell }(x)\ts\mathcal{A}_{\ell }(y)} 
\K (u,v) \mmm (du) \mmm (dv)
,\end{align}
where $ \mathcal{A}_{\ell }(x) $ is such that 
$ x \in \mathcal{A} _{\ell }(x) \in \Delta (\ell )$. Then 
$ \K _{\ell } $ can be regarded as a discrete kernel on 
$ \Il $. If $ \K _{\ell } $ satisfies \thetag{A1}, then $ \K _{\ell } $ generates deteminantal point field $ \mu_{\K _{\ell } }$. 
Indeed, $ \mathrm{Spec}(\K _{\ell } ) \subset [0,1]$ follows from 
$ \mathrm{Spec}(\K ) \subset [0,1]$ and the Fubini theorem. 
One can expect the convergence of the kernel $ \K _{\ell } $ 
to $ \K $, and as a result, the weak convergence of $ \mu_{\K _{\ell } }$ 
to $ \mu $, at least for continuous $ \K $. 
Because $ \mu_{\K _{\ell } }$ is a determinantal point process on the discrete space, 
its tail $ \sigma $-field is trivial. 
Such weak convergence, however, does not suffice for 
the convergence of the values on the tail $ \sigma $-field $ \TailS $. 

Taking the above into account, we consider the second approximation 
given by $ \muls $ below.  
Let $ \Gl $ be the sub-$ \sigma $-field of $ \mathcal{B} (\mathsf{S})$ given by 
\begin{align}\label{:20e}& 
\Gl  =  \sigma [ \{\mathsf{s}\in \mathsf{S}; \mathsf{s}(\Ali ) = n 
\}; i\in \Il , n \in \mathbb{N} ]
.\end{align}
Combining \eqref{:10b} and \eqref{:10c} with \eqref{:20e}, we obtain 
\begin{align}\label{:20f}& 
\Gl \subset \Gll ,\ \quad 
\sigma [\Gl; \ell \in \N ] = \mathcal{B}(\Ssf ) 
.\end{align}
Let $ \muls  $ be the regular conditional probability of 
$ \muk $ with respect to $ \Gl  $. 
Using \eqref{:20f}, we shall prove in \lref{l:61} that 
for all $ U \in \mathcal{B}(\Ssf ) $
\begin{align}\label{:11b}&
\limi{\ell } \mulsU = 1_{ U } (\mathsf{s}) 
\quad \text{ for $ \muk $-a.s. } \mathsf{s}
.\end{align}
We see that the convergence in \eqref{:11b} is stronger than the weak convergence. 
In particular, the convergence in \eqref{:11b} is valid for all 
$ U \in \TailS $ because $ \TailS  \subset \mathcal{B}(\Ssf ) $. 

We can naturally regard $ \Dl = \{ \Ali  \}_{i \in \Il } $
 as a discrete, countable set 
with the interpretation that each element $ \Ali  $ is a point. 
Thus, $ \muls $ can be regarded as a point process 
on the discrete set $ \Dl $.

If $ \muls $ were a determinantal point process for each $ \ell $, 
then \tref{l:11} would follow from \eqref{:11b} immediately 
because determinantal point processes on discrete spaces 
always have trivial tails, and as discussed above, 
$ \muls $ is naturally regarded as a determinantal point process 
on the discrete space $  \Dl $. This is clearly not the case because 
determinantal point processes are supported on single configurations and 
\begin{align}\label{:11d}&
 \mu  (\{ \sss ; \sss (\Ali ) \ge 2 \}|\Gl ) > 0 
.\end{align}
Hence we introduce a sequence of 
{\it fiber bundle-like sets} $ \IIIl $ ($\ell \in \N )$ in \sref{s:2} 
with base space $ \Dl $ with fiber consisting of a set of binary trees. 
We further expand $ \IIIl $ to $ \Omega (\ell )$ in \eqref{:22t}, which has a fiber 
whose element is a product of a tree $ \ii $ and a component 
$ \Bli $ of partitions. See notation after \tref{l:21}. 

Let $ \mul $ denote the restriction of $ \mu $ on $ \Gl $. 
By construction $ \mul (\mathsf{A}) = \mulA $ for all $ \mathsf{A} \in \Gl $. 
In \tref{l:21} and \tref{l:22}, we construct a lift $ \7 $ of $ \mul $ 
on the fiber bundle $ \Omega (\ell )$, and prove  tail triviality of 
the lift $ \7 $ in \tref{l:24}, 
which establishes  tail triviality of $ \mul $ in \tref{l:25}. 
Combining \tref{l:25} with the martingale convergence theorem in \lref{l:61}, 
we obtain \tref{l:11}.  

The key point of the construction of the lift $ \7 $ is that we construct 
a consistent family of orthonormal bases $ \6  $ 
in \eqref{:21j} and \eqref{:21k}, and that we introduce the kernel 
 $ \KFl $ on $ \IIIl $ in \eqref{:21v} such that 
\begin{align}\tag{\ref{:21v}}&
\KFl (\ii ,\jj ) = 
\int_{\Sit \ts \Sit } \K (x,y) 
\overline{\flii  (x) } \fljj  (y) \mmm (dx)\mmm (dy) 
.\end{align}
We shall prove in \lref{l:32} that $ \KFl $ is a determinantal kernel on 
$ \IIIl $, 
and present $ \nuFl $ as the associated determinantal point process on 
$ \IIIl $. 
To some extent, $ \nuFl $ is isometric to $ \mul $ 
through the orthonormal basis $ \6  $. 
We shall indeed prove in \tref{l:21} that their correlation functions 
$ \rhol ^m $ and $ \rhoFl ^m $ satisfy the identity: 
\begin{align}\tag{\ref{:21b}}&
\int_{\mathbb{A}} \rhol ^m(\mathbf{x}) \mmm ^{m}(d\mathbf{x}) = 
\sum_{\iii \in \III_{\ell }(\mathbb{A}) } 
\rhoFl ^m (\iii ) 
,\end{align}
which is a key to construct the lift $ \7 $.

While preparing the manuscript, we have heard that Professor A. Bufetov 
has proved independently tail triviality of determinantal point processes on continuous spaces independently of us 
(a seminar talk at Kyushu University in October 2015). 
His method is completely different from ours and requires a restriction on 
an integrability condition of the determinantal kernel $ \K (x,y)$. 
An improved version of the work is now available in \cite{bqs}.

The organization of the paper is as follows. 
In \sref{s:2}, we introduce definitions and concepts and state the main theorems (Theorems \ref{l:21}--\ref{l:25}). We give tree representations of $ \mu $. 
In \sref{s:3}, we prove \tref{l:21}. 
In \sref{s:4}, we prove \tref{l:22}--\tref{l:25}. 
In \sref{s:6}, we prove \tref{l:11}. 
In \sref{s:7}, we present motivational examples such as Sine$_2 $, Airy$_2 $, and Bessel$_2 $, and Ginibre point processes.

\section{Set up and main results}\label{s:2}

In this section, we recall various essentials and present the main theorems 
(\tref{l:21}--\tref{l:25}) other than \tref{l:11} . 

A configuration space $ \Ssf $ over $ \Sit $ 
is a set consisting of configurations on $ \Sit $ such that 
\begin{align}& \label{:20a}
\Ssf = \{ \mathsf{s}\, ;\,   \mathsf{s}=\sum_i \delta _{s_i},\, 
\{ s_i \}\subset \Sit ,\, 
\mathsf{s}(K) < \infty \text{ for any compact } K \} 
,\end{align}
where $ \delta_{s_i}$ denotes the delta measure at $ s_i$. 
A probability measure $ \mu $ on $ ( \Ssf , \mathcal{B}(\Ssf ) )$ 
is called a point process, also called random point field. 
 A symmetric function $\rho ^m $ on $\Sit ^m$ is called 
the $ m $-point correlation function of a point process $ \mu $ 
with respect to a Radon measure $ \mmm $ if it satisfies 
\begin{align} \label{:20b}
\int_{\Ssf } 
\prod_{i=1}^j \frac{\mathsf{s}(A_i)!}{(\mathsf{s}(A_i) - k_i)!} 
\mu(d\mathsf{s}) 
&= \int_{A_1^{k_1} \times \cdots \times A_j^{k_j}} 
\rho ^m (\mathbf{x}) \mmm ^m(d\mathbf{x})
.\end{align}
Here $A_1, \dots, A_j \in \mathcal{B}(\Sit )$  are disjoint and 
$k_1,\dots, k_j \in \N $ such that $k_1+\cdots + k_j = m $. 
If $\mathsf{s}(A_i) - k_i \le 0$, we set 
${\mathsf{s}(A_i)!}/{(\mathsf{s}(A_i) -k_i)!}=0$.

We fix a point $ o \in \Sit $ as the origin, and set 
$ S_r = \{ x\in \Sit \, ;\, \dist (o,x) < r \} $. 
Each $ S_r $ is assumed to be relatively compact, and thus 
$ \sss (S_r ) < \infty $ for all $ \sss \in \Ssf $ and $ r \in \N $. 
In this sense, each element $ \sss $ of $ \Ssf $ is 
a locally finite configuration. We note that this notion depends on the choice of metric $ \dist $ on $ \Sit $. 

For a Borel set $ A $ we set $ \map{\pi_{A}}{\Ssf }{\Ssf }$ 
by $ \pi_{A}(\mathsf{s}) (\cdot ) = \mathsf{s} (\cdot \cap A )$. 
We set $ \map{\pi_{S_r^c}}{\Ssf }{\Ssf }$ 
such that $ \pi_{S_r^c} (\mathsf{s}) = \mathsf{s}(\cdot \cap S_r^c )$. 
We denote by $ \TailS  $ the tail $ \sigma $-field such that 
\begin{align}\label{:20c}&
\TailS  = \bigcap_{\rr = 1}^{\infty} \sigma [\pi_{S_r^c}]
.\end{align}
If we replace $ S_r $ by any increasing sequence $ \{O_{\rr  }\}$ of relatively compact open sets such that 
$ \cup_{\rr  = 1}^{\infty} O_{\rr  } = \Sit $, 
then $ \TailS  $ defines the same $ \sigma $-field. 
Thus $ \TailS  $ is independent of the choice of $ \{O_{\rr  }\}$. 

Let $ \Dl =\{ \Ali  \}_{i\in \Il  } $ be as in \thetag{A2}, where 
 $ \ell \in \N $. 
We set $ \Delta =  \{ \mathcal{A}_{\ii }  \}_{\ii \in I } $ such that 
 \begin{align}\notag &
 \Delta  = \Done , \quad \mathcal{A}_{\ii }=\mathcal{A}_{1,\ii  }  \quad I = I (1)
 .\end{align}
In consequence of  \eqref{:10d}, we assume 
without loss of generality that each element $ \ii  $ of 
the parameter set $ \Il $ is of the form 
\begin{align}\label{:20g}&
\Il = \I \ts \{ 0,1 \}^{\ell -1}
.\end{align}
That is, each $ \ii \in \Il $ is of the form 
$ \ii  = (j_1,\ldots,j_{\ell }) \in \I \ts \{ 0,1 \}^{\ell -1} $. 
We take a label $ \ii \in \cup_{\ell = 1}^{\infty} \Il $ in such a way that, 
for $ \ell < \ell ' $, $ \ii \in \Il $, and $ \ii ' \in I (\ell ' )$, 
\begin{align}\notag &
\Ali  \supset \AAAA _{\ell ', \ii '}  
\Leftrightarrow 
\ii =(j_1,\ldots,j_{\ell }) \text{ and }
\ii '=(j_1,\ldots,j_{\ell },\ldots,j_{\ell '})
.\end{align}
We denote by $ \wIII $ the set of all such parameters:
\begin{align}\label{:20i}&
\wIII  = \bigcupii{\ell } \Il = \bigcupii{\ell }  \I  \ts \{ 0,1 \}^{\ell -1} 
.\end{align}
We can regard $ \wIII $ as 
a collection of binary trees and $ \I  $ is the set of their roots. 

\begin{example}[Binary partitions of $ \mathbb{R}$] \label{d:21}
Typically we can take $ \Sit = \R $, $ \mmm (dx)= dx $, and $ \I  = \mathbb{Z}$. 
For $ \ii = (j_1,\ldots,j_{\ell }) \in \Il $, we set 
$ \Jonei = j_1 $ and, for $ \ell \ge 2 $, 
\begin{align}\label{:20q}&
\Jli  = j_1+\sum_{n=1}^{\ell - 1} \frac{j_n}{2^{n}}
.\end{align}
We take $ \Ali = [\Jli , \Jli  + 2^{-\ell +1})$. 
\end{example}

For $ \ii  = (j_1,\ldots,j_{\ell }) \in \wIII $, we set 
$ \mathrm{rank}(\ii ) = \ell $. 
For $ \ii $ with $ \mathrm{rank}(\ii ) = \ell $, 
we set 
\begin{align}\label{:20j} &\Bi  =
\begin{cases}
 \mathcal{A}_{1,\ii } & \ell = 1,\\
\mathcal{A}_{\ell -1, \ii ^-}  
&  \ell \ge 2
,\end{cases}
\end{align}
where $ \ii ^- = (j_1,\ldots,j_{\ell -1})$ for 
$ \ii = (j_1,\ldots, j_{\ell })\in \Il $. 
Let $ \III \subset \wIII $ such that 
\begin{align}\label{:20k}&
\III = \I \cup \Big\{  \bigcup_{\ell = 2}^{\infty}
\{ \ii  \in \Il ; j_{\ell }=0   \}\Big\}
,\end{align}
where $ \ii = (j_1,\ldots,j_{\ell }) \in \Il $. 

Let $ \mathbb{F} = \{ \fii \}_{\ii \in\mathbb{I} } $ 
be an orthonormal basis of $ L^2(\Sit ,\mmm ) $ satisfying 
\begin{align}
& && \label{:20l}
\sigma [\fii ; \ii \in \III ,\, \mathrm{rank} (\ii ) = \ell ] = \Fl 
&&\text{ for each $ \ell \in \N $} 
,\\ \label{:20m}&&&
\mathrm{supp}(\fii ) 
= 
\Bi  
&& \text{ for each $ \ii \in \III $}
,\\
\label{:20n}&&&
\fii  (x) = 
1_{\Ai }(x)/\sqrt{\mmm (\Ai )}
&&\text{ for $ \mathrm{rank} (\ii ) = 1$}
.\end{align}
For a given sequence of $ \mmm $-partitions 
satisfying \thetag{A2}, such an orthonormal basis exists. 
We present here an example. 

\begin{example}[Haar functions] \label{d:22} 
We make the same assumptions as in \eref{d:21}. 
Let $ \ii = (j_1,\ldots,j_{\ell }) \in \III $. We set for, $ \ell = 1 $ and $ \ii = (j_1)$, 
\begin{align}\notag & 
f_{\ii } (x) = 
1_{[j_1,j_1+1)}  (x) 
\end{align}
and, for $ \ell \ge 2 $ and $  \ii = (j_1,\ldots,j_{\ell }) \in \III $, 
\begin{align}\notag & 
f_{\ii } (x) = 
2^{(\ell -1)/2 } \{ 1_{[\Jli , \Jli  + 2^{-\ell +1})}(x) 
- 
1_{[\Jli + 2^{-\ell +1}, \Jli  + 2^{-\ell +2 })}(x) \}  
.\end{align}
We can easily see that $ \{ \fii  \}_{\ii \in \III } $ is an orthonormal basis of 
$ L^2(\R , dx )$. 
We remark that $ j_{\ell } = 0 $ because $  \ii = (j_1,\ldots,j_{\ell }) \in \III $ as we set in \eqref{:20k}.  
\end{example}

We next introduce the $ \ell $-shift of above objects 
such as $ \III $, $ \Bi $, and $ \mathbb{F} = \{ \fii \}_{\ii \in\mathbb{I} } $. 
Let $ \wIII (1) = \wIII $ and, for $ \ell \ge 2 $, 
\begin{align}\label{:21g}
\quad \quad \wIIIl := &
 \bigcupii{\rfrak } 
\Il \ts \{ 0,1 \}^{\rfrak -1 } 
,\end{align}
where $ \Il =  \I \ts \{ 0,1 \}^{\ell -1}$ is as in \eqref{:20g}.  
For $ \ell , \rfrak \in  \N  $, we set 
$ \map{\theta _{\ell -1 ,\rfrak }}{\wIII }{ \wIIIl }$ 
such that 
$ \theta _{0,\rfrak } = \mathrm{id}$ ($ \ell = 1$) and, for $ \ell \ge 2$,   
\begin{align}\label{:21h}&
\theta_{\ell -1,\rfrak } ((j_1,\ldots,j_{\ell + \rfrak -1 })) 
= ( \mathbf{j}_{\ell }, j_{\ell +1 }, \ldots,j_{\ell + \rfrak -1 }) 
\in \Il  \ts \{ 0,1 \}^{\rfrak -1 } 
,\end{align}
where $ \mathbf{j}_{\ell }=(j_1,\ldots,j_{\ell }) \in \Il $. 
For $ \ell = 1$, we set $ \III (1) = \III $. For $ \ell \ge 2 $, we set 
\begin{align}\label{:21i}&
\IIIl  = \Il \cup \Big\{  
\bigcup_{\rfrak = 2}^{\infty}  \theta_{\ell -1 ,\rfrak }  (\III )
\Big\}
.\end{align}

We set $ \mathrm{rank}(\ii ) = r $ for $ \ii \in  \Il  \ts \{ 0,1 \}^{\rfrak -1 } $. 
By construction $ \mathrm{rank}(\ii ) = r $ 
for $ \ii \in \theta_{\ell -1,\rfrak } (\wIII )$. 
Let $ \6  $ such that, 
for $ \rfrak = \mathrm{rank} (\ii ) $,  
\begin{align}\label{:21j}&&&\quad \quad 
\flii  (x) = 
1_{\Ali }(x)/\sqrt{\mmm (\Ali )}
&&\text{ for } \rfrak = 1
,\\
&&& \label{:21k} \quad \quad  \flii (x) = 
f_{ \theta_{\ell -1 ,\rfrak } ^{-1}(i)} (x) 
&&
\text{ for } \rfrak  \ge 2 
,\end{align}
where $ \Dl = \{ \Ali  \}_{i\in\Il } $ is given in \thetag{A2}. 
Then $ \6  $ is an orthonormal basis of $ L^2(\Sit ,\mmm ) $. 
This follows from assumptions \eqref{:21j} and \eqref{:21k} and 
the fact that $ \FFF = \{ \fii \}_{\ii \in\mathbb{I} } $ 
is an orthonormal basis.

\begin{remark}\label{r:21} \thetag{1} 
We note that 
$ \flii \in \FFFl $ is a newly defined function 
if $ \mathrm{rank}(\ii ) = 1 $, whereas 
$ \flii \in \FFFl $ is an element of $ \FFF $ if 
$ \mathrm{rank}(\ii ) \ge 2 $. 
In particular, we see that 
\begin{align}\label{:21l}&
\{ \flii  \}_{\ii \in \IIIl , \, \mathrm{rank}(\ii )\ge 2} 
\subset 
\{ \fii  \}_{\ii \in \III , \, \mathrm{rank}(\ii )\ge 2} 
.\end{align}
\thetag{2} 
Let 
$ \jj = (j_1,\ldots,j_{\ell + \rfrak -1})\in\III $ 
and 
$ \ii  = (\mathbf{j}_{\ell }, j_{ \ell + 1 },\ldots,j_{\ell + \rfrak -1}) \in \IIIl $. 
Then 
$$ \jj =  \theta_{\ell -1 ,\rfrak } ^{-1}(\ii )
.$$
Furthermore, $ \flii  \in\FFFl  $ and $ \fjj \in  \FFF $ satisfy $ \flii = \fjj   $ 
for $ r=\mathrm{rank} (\ii ) \ge 2 $. 
\\
\thetag{3} By construction, we see that 
\begin{align} 
& && \label{:21s}
\sigma [\flii ; \ \ii \in \IIIl , \, \mathrm{rank} (\ii ) = r ] = 
\mathcal{F}_{\ell -1 + r}  
&&\text{ for each $ \ell , r \in \N $} 
,\\ \label{:21t}&&&
\mathrm{supp}(\flii ) 
 = \Bli 
&&
\text{ for all } \ii \in \IIIl 
,\end{align}
where we set, for $  \jj = \theta_{\ell -1 ,\rfrak } ^{-1}(\ii )$ 
such that $ \mathrm{rank}(\ii ) = r $, 
\begin{align}\label{:21u}&
\Bli  = \Bj  
.\end{align}
\end{remark}

Using the orthonormal basis 
$ \6  $, 
we set $ \KFl $ on $ \IIIl $ by 
\begin{align}\label{:21v}&
\KFl (\ii ,\jj ) = 
\int_{\Sit \ts \Sit } \K (x,y) 
\overline{\flii  (x) } \fljj  (y) \mmm (dx)\mmm (dy) 
.\end{align}
Let $ \lambdal $ be the counting measure on $ \IIIl $. 
We shall prove in \lref{l:32} that $ (\KFl , \lambdal )$ satisfies \thetag{A1}. 
Hence we obtain the associated determinantal point process $ \nuFl $ on $ \IIIl $ 
from general theory \cite{sosh,s-t.aop}.

For $ \iiM  \in \IIIl $, let $ \mmFi  (dx ) $ be 
the probability measure on $ \Sit $ such that 
\begin{align}\label{:21w}&
\mmFi  (dx) =  | \flii (x)|^2 \mmm (dx) 
.\end{align}
For $ \iii = (\ii _n)_{n=1}^m \in \IIIl ^m $ and 
$ \mathbf{x} = (x_n)_{n=1}^m  $, where $ m \in \N \cup \{ \infty  \} $, we set 
\begin{align}\label{:21x}&
\mmFiii  (d\mathbf{x}) = \prod_{n=1}^m 
 |f_{\ell ,\ii  _n }(x_n)|^2 \mmm (dx_n) 
.\end{align}
By \eqref{:21k} $ \mmFiii $ is a probability measure on $ \Sit ^m $. 
By \eqref{:21t}, we have 
\begin{align}& \label{:21y}
\mmFiii   ( \9 ) = 1 
.\end{align}

Let $ \Gl $ be the sub-$ \sigma $-field  as in \eqref{:20e}.  
Let $ \nul $ be the $ (\KFl , \lambdal )$-determinantal point process 
as before. 
Let $ \rhol ^m $ and $ \rhoFl ^m $ be 
the $ m $-point correlation functions of $ \mul $ and 
$ \nuFl $ with respect to $ \mmm $ and $ \lambdal $, respectively. 
We now state one of our main theorems:

\begin{theorem} \label{l:21} 
Let $ \III_{\ell }(\AAAA ) = \{ \ii \in \IIIl \, ;\,  \Bli \subset \AAAA \} $. 
For $ \mathbb{A}=  \AAAA _1\ts \cdots \ts \AAAA _m$, we set 
\begin{align}\label{:21a}&
\III_{\ell }(\mathbb{A}) = 
\III_{\ell }(\AAAA _1)\ts \cdots \ts \III_{\ell }(\AAAA _m ) 
.\end{align}
Assume that 
$\AAAA _n \in \Dl  $ for all $ n = 1,\ldots,m $. 
Then
\begin{align}\label{:21b}&
\int_{\mathbb{A}} \rhol ^m(\mathbf{x}) \mmm ^{m}(d\mathbf{x}) = 
 \sum_{\iii \in \III_{\ell }(\mathbb{A}) } 
\rhoFl ^m (\iii ) 
.\end{align}
\end{theorem}

Let $ \mathsf{I}(\ell ) $ be the single configuration space over $ \IIIl $. 
We write $ i \in \mathsf{i}$ if $ \mathsf{i}(\{ i \} ) = 1 $. 
Each 
$ \mathsf{i}=\sum_{i\in\mathsf{i}} \delta_{i}\in  \mathsf{I}(\ell ) $ 
can be regarded as a subset of $ \IIIl $ 
by the correspondence of $ \mathsf{i} $ to $ \{ i \}_{i\in\mathsf{i}} $. 
Let 
\begin{align}& \label{:22t}
 \1 (\ell )  := 
\bigcup_{i \in \IIIl } \{ i \} \ts \Bli  
.\end{align}
Let $ \4   $ be the single configuration space over $ \1 (\ell )$. 
Then by definition each element $ \omega \in \4   $ is of the form 
$ \omega = \sum_{i\in\mathsf{i}} \delta_{(\ii , s_{\ii })}$ 
such that $ s_{\ii } \in \Bli $. Hence 
\begin{align} \label{:22u}
\4 &  \subset \{ \omega = \sum_{i\in\mathsf{i}} \delta_{(\ii , s_{\ii })} 
\, ;\, \mathsf{i} =\sum_{i\in\mathsf{i}} \delta_{i}\in \IIl , \, 
 s_{\ii } \in \Bli     \} 
.\end{align}
Let $ \mmFi $ be as in \eqref{:21w}. We set  
\begin{align}
\label{:22v}
\mmFl = & 
 \prod_{{ \ii \in \IIIl  }}  \mmFi  
,\quad 
\mmFii = 
 \prod_{{ \ii \in \mathsf{i}  }}  \mmFi  
.\end{align}

\begin{remark}\label{r:22}
Let $ \mathbf{i}= (i_1,\ldots,i_m) $ and 
$ \mathsf{i}= \sum_{n=1}^m \delta_{i_n}
\equiv \sum_{i\in\mathsf{i}} \delta_{i}$. 
By definition $ \mmFii $ in \eqref{:22v} 
is a product measure on the product space 
$ \prod_{\ii \in \mathsf{i}} \Bli $ 
with (unordered) parameter $ \ii \in \mathsf{i}$, whereas 
$ \mmFiii $ in \eqref{:21x} is a product measure on the product space 
$ \B _{\ii _1}\ts \cdots \ts \B _{\ii _m} $ 
with (ordered) parameter $ \mathbf{i}= (i_1,\ldots,i_m) $. 
\end{remark}

We set $ \map{ \iotal  }{ \4 }{ \mathsf{I}(\ell ) }$ 
such that $ \iotal  ( \omega ) = \mathsf{i}$, and 
$ \map{ \kappal }{ \4 }{  \prod_{\ii \in \mathsf{i}} \Bli }$ 
such that $ \kappal (\omega ) = (s_i)_{i\in\mathsf{i}}$, 
where  
$ \omega =  \sum_{i\in\mathsf{i}}  \delta_{(\ii , s_i )} $, 
$ \mathsf{i} =  \sum_{i\in\mathsf{i}} \delta_{i} $, and 
$ s_i \in \Bli $.

Let $ \7 $ be the probability measure  on $ \4 $ given by   
\begin{align} \label{:22w} & \quad \quad 
(\7 ) \circ \iotal ^{-1} (d\mathsf{i}) = \nul (d\mathsf{i})
,\\\label{:22x}&\quad \quad 
\7 ( \kappal (\omega ) \in d\mathbf{s}  | \iotal (\omega ) = \mathsf{i}) = 
   \mmFii (d\mathbf{s})  
,\quad \mathbf{s}=(s_i)_{i\in\mathsf{i}}
.\end{align}

\begin{remark}\label{r:23} \thetag{1} 
We can naturally regard the probability measures in \eqref{:22x} 
as a point process on $ \prod_{\ii \in \mathsf{i}} \Bli $ 
supported on the set of configurations with exactly one particle configuration 
$ \mathsf{s}= \delta_{\mathbf{s}}$ on 
$ \prod_{\ii \in \mathsf{i}} \Bli $, that is, 
$ \mathbf{s}=(s_i)_{i\in\mathsf{i}}$ is such that  
$ s_i \in \Bli $ for each $ i\in\mathsf{i}$. 
\\\thetag{2} 
We can regard $ \7 $ as a marked point process as follows: The configuration $ \mathsf{i}$ is distributed according to $ \nul $, while the marks are independent and for each $ \mathsf{i}$ the mark $ \mathbf{s}$ 
is distributed according to $ \mmFii $. Thus the space of marks depends on $ \mathsf{i}$. 

\end{remark}

\begin{theorem}	\label{l:22} 
Let $ \map{\ulabl  }{\4 }{\Ssf }$ be such that 
$ \ulabl (\omega ) = \sum_{\ii \in\mathsf{i}} \delta_{s_{\ii }}$, 
where $ \omega = \sum_{\ii \in\mathsf{i}} \delta_{(\ii , s_{\ii } )} $. 
Then 
\begin{align}\label{:22a}&
\mul  = (\7 )\circ  \ulabl  ^{-1} |_{ \Gl } 
.\end{align}
\end{theorem}

\begin{remark}\label{r:4} 
\tref{l:22} implies that $ \7 $ is a {\em lift} of $\mul $ onto $ \4 $. 
We can naturally regard  $ \wIIIl  $ as binary trees. 
Hence we call $ \7 $ a tree representation of $ \mu $ of level $ \ell $. 
\end{remark}

We present a decomposition of $ \mul $, 
which follows from \tref{l:22} immediately. 
Let $ \mmFii ^{\ulab } = \mmFii \circ \ulabi ^{-1}$, where 
$ \map{\ulabi }{ \prod_{\ii \in \mathsf{i}} \Bli  }{\Ssf }$ 
is the unlabel map such that 
$ \ulabi ((s_{\ii })_{\ii \in \mathsf{i}}) = 
\sum_{\ii \in \mathsf{i}} \delta_{s_{\ii }}$. 
\begin{theorem}	\label{l:23} 
For each $ \mathsf{A}\in \Gl $, 
\begin{align}\label{:23b}& \quad \quad 
\mu (\mathsf{A}) = 
\int_{\IIl }\nuFl (d\mathsf{i}) \,  \mmFii ^{\ulab }    (\mathsf{A}) 
.\end{align}
\end{theorem}

Let 
$ \IIIl _{p} = \{ \ii \in \IIIl ; 
r \le p , |j_1|\le p \} $, where $ i =(\mathbf{j}_{\ell },j_{\ell+1},\ldots,j_{\ell+r-1})$, $ r = \mathrm{rank}(\ii )$, and 
$ \mathbf{j}_{\ell }= (j_1,j_2,\ldots,j_{\ell })$. 
Let 
$ \pi_{p}^c (\mathsf{i}) = \mathsf{i} (\cdot \cap \IIIl _{p}^c)$. 
Then we set $ \Tail (\IIIl ) = \cap_{p=1}^{\infty}\sigma [\pi_{p}^c ]$. 
From this we can define the tail $ \sigma $-field 
$ \Tail (\4  )$ of $ \4  $ 
because $ \Omega (\ell ) $ is a subset of $ \IIIl \ts \Sit $. 
\begin{theorem} \label{l:24} 
$ \7 $ is trivial on 
$ \Tail (\4 ) \cap \ulabl ^{-1} (\Gl )$. That is, 
\begin{align}\label{:24a}& \quad \quad \quad 
\7 (\mathsf{A}) \in \{ 0,1 \} 
\quad \text{ for all $ \mathsf{A} \in 
 \Tail (\4 ) \cap \ulabl ^{-1} (\Gl ) $}
.\end{align}
\end{theorem}

We remark that $ \mul $ is {\em not} a determinantal point process. 
Hence we exploit $ \7 $ instead of $ \mul $. 
As we have seen in \tref{l:22},  
$  \7 $ is a lift of $ \mul $ in the sense of \eqref{:22a}, 
from which we can deduce nice properties of $ \mul $. 
Indeed, an application of \tref{l:22} combined with \tref{l:24} 
is  tail triviality of $ \mul $: 
\begin{theorem}	\label{l:25}
$ \mul $ is tail trivial. That is, 
\begin{align}\label{:25a}&
\mul (\mathsf{B}) \in \{ 0,1 \} 
\quad \text{ for all $ \mathsf{B} \in \TailS \cap \Gl $}
.\end{align}
\end{theorem}

We shall apply \tref{l:25} to prove \tref{l:11} in \sref{s:6}.

\section{Proof of \tref{l:21}}\label{s:3}
The purpose of this section is to prove \tref{l:21}. 
In \lref{l:31}, we present 
the identity of kernels $ \K $ and $ \KFl $ 
using the orthonormal basis $ \FFFl  $, 
where $ \KFl $ is the kernel given by \eqref{:21v} and 
$ \FFFl  $ is as in \eqref{:21j} and \eqref{:21k}. 
In \lref{l:32}, we prove $ (\KFl , \lambdal )$ is a determinantal kernel and 
the associated determinantal point process $ \nul $ exists. 
We will lift the the identity between $ \K $ and $ \KFl $ 
to that of correlation functions of $ \mul $ and $ \nul $ in \tref{l:21}.  

%

By definition $ \6 $ satisfies 
\begin{align}\label{:30a}&
\int_{\Sit } |\flii (x)|^2 \mmm (dx)= 1 \quad 
\text{ for all } \ii \in\IIIl 
,\\\label{:30b}&
\int_{\Sit } \flii (x) \overline{\fljj (x)}  \mmm (dx)= 0 \quad 
\text{ for all } \ii \not=\jj \in \IIIl 
.\end{align}
\begin{lemma} \label{l:31} 
\thetag{1} Let $ P(x) = \sum_i \xi (\ii ) f _{\ell ,i} (x) $ and 
$ Q(y) = \sum_j \eta (\jj ) f _{\ell ,j} (y) $. 
Suppose that the supports of $ \xi $ and $ \eta $ are finite sets. 
Then 
\begin{align}\label{:31z}&
\int_{\Sit \ts \Sit }   \K (x,y) \overline{P (x)}Q (y) \mmm (dx) \mmm (dy)  = 
\sum_{i,j} \KFl (i,j) \overline{\xi (\ii )}\eta (\jj ) 
.\end{align}
\thetag{2} We have an expansion of 
$ \K $ in $ L^2_{\mathrm{loc}}(\Sit \ts \Sit ,\mmm \ts \mmm )$ such that 
\begin{align}\label{:31a}&
\K (x,y) = \sum_{\ii ,\jj \in\IIIl }
\KFl (\ii ,\jj ) 
\flii  (x) 
\overline{\fljj  (y)} 
.\end{align}
\end{lemma}
\proofbegin 
From \eqref{:21v} we deduce that  
\begin{align}\label{:31b}&
\int_{\Sit \ts \Sit }   \K (x,y) 
\overline{P (x)}Q (y) \mmm (dx) \mmm (dy) 
\\ \notag =& 
\int_{\Sit \ts \Sit }  \K (x,y) 
\overline{\sum_i \xi (\ii ) f _{\ell ,i} (x) } 
\sum_j \eta (\jj ) f _{\ell ,j} (x) \mmm (dx) \mmm (dy) 
\\ \notag =&
\sum_{i,j} \int_{\Sit \ts \Sit } 
 \K (x,y)  \overline{\flii (x)}  \fljj (y) \mmm (dx) \mmm (dy) 
\overline{\xi (\ii )}\eta (\jj ) 
\\ \notag =&
\sum_{i,j} \KFl (i,j) \overline{\xi (\ii )}\eta (\jj ) 
.\end{align}
This yields \eqref{:31z}. We have thus proved \thetag{1}. 
By a direct calculation, we have 
\begin{align}\label{:31c}&
 \int_{\Sit }  P(x) \overline{\flii (x)}  \mmm (dx) = 
\int_{\Sit }  
\sum_p \xi (p) f _{\ell ,p} (x) \overline{\flii (x)}  \mmm (dx) = 
\xi (\ii ) 
,\\ \notag &
 \int_{\Sit }  Q(y) \overline{\fljj (y)}  \mmm (dy) = 
\int_{\Sit }  \sum_q \eta (q) f _{\ell ,q} (y) \overline{\fljj (y)}  \mmm (dy) = \eta (\jj ) 
.\end{align}
Combining \eqref{:31b} and \eqref{:31c} yields 
\begin{align}\notag 
&
\int_{\Sit \ts \Sit }   \K (x,y) \overline{P (x)}Q (y) \mmm (dx) \mmm (dy) = 
\\ \notag &
\int_{\Sit \ts \Sit }  \sum_{i,j} \KFl (i,j) 
\flii (x) \overline{\fljj (y)} \overline{P (x)}Q (y) \mmm (dx) \mmm (dy) 
.\end{align}
This implies \eqref{:31a}.  
\proofend

Let $ \lambdal $ be the counting measure on $ \IIIl $ as before. 
We can regard $ \KFl $ as an operator on 
$ L^2(\IIIl  ,\lambdal ) $
such that $ \KFl \xi (\ii ) = \sum_{\jj \in \IIIl } \KFl (\ii ,\jj ) \xi (\jj ) 
$. 
We now prove that the $ (\KFl ,\lambdal ) $-determinantal point $ \nul $
process exists. 
\begin{lemma} \label{l:32} 
Let $ \mathrm{Spec}( \KFl  ) $ be the spectrum of $ \KFl $. Then  
\begin{align}\label{:32a}&
 \mathrm{Spec}( \KFl  ) \subset [0,1]
.\end{align}
In particular, there exists a unique, determinantal point process $ \nuFl $ 
on $ \IIIl $ associated with $ (\KFl ,\lambdal ) $. 
\end{lemma}
\proofbegin 
Recall that $\6 $ is an orthonormal basis of $ L^2(\Sit ,\mmm )$. 
Let $ \map{U}{L^2(\Sit ,\mmm )}{L^2(\IIIl , \lambdal  )}$ be the unitary operator such that 
$ U (\flii ) = e_{\ell ,i}$, where 
$ \{e_{\ell ,i} \}_{i\in\IIIl }$ 
is the canonical orthonormal basis of $L^2(\IIIl , \lambdal  )$. 
Then by \lref{l:31} we see that $  \KFl  = U \K U^{-1}$. 
Hence $ \KFl $ and $ \K $ have the same spectrum. 
We thus obtain \eqref{:32a}.  
Because $ \KFl  $ is Hermitian symmetric, 
the second claim is clear from \eqref{:32a}, \thetag{A1}, 
and general theory 
\cite{sosh,s-t.jfa,s-t.aop}. 
\proofend

\begin{lemma} \label{l:33} 
Let $ \Bli = \mathrm{supp} (\flii )$ be as in \eqref{:21t}.  
Then, for $ \ii , \jj  \in \IIIl $ and $ \AAAA \in \Fl  $, 
\begin{align}\label{:33a}&
\int_{ \AAAA } \flii (x)\overline{\fljj (x)} \mmm (dx) =
\begin{cases} 1 & (\ii = \jj ,\ \Bli  \subset \AAAA )
\\
0 & (\text{otherwise})
\end{cases}
.\end{align}
\end{lemma}
\proofbegin  
We recall that $\Bli $ is the support of $ \flii $ by \eqref{:21t}. 
Suppose $ \ii = \jj $ and $ \Bli \subset \AAAA $. Then 
from \eqref{:30a} 
\begin{align}\label{:33b}&
\int_{ \AAAA } \flii (x)\overline{\fljj (x)} \mmm (dx) = 
\int_{\Sit } \flii (x)\overline{\flii (x)} \mmm (dx) = 1 
.\end{align}
Suppose that $ \ii = \jj $ and that $ \Bli \not\subset \AAAA $. 
Then, using $ \AAAA \in \Fl  $, \eqref{:20j}, and \eqref{:21u}, 
we deduce that $ \Bli \cap \AAAA = \emptyset $. 
Because $ \Bli = \mathrm{supp} (\flii ) $, we obtain 
\begin{align}\label{:33c}&
\int_{ \AAAA } \flii (x)\overline{\fljj (x)} \mmm (dx) = 0 
.\end{align} 
Finally, suppose $ \ii \not= \jj $. Because $ \AAAA \in \Fl  $,  
we see that $ \Bli \subset \AAAA $ or $ \Bli \cap \AAAA = \emptyset $. 
The same also holds for $ \Blj $. In any case, we obtain \eqref{:33c} 
from \eqref{:30b}.  
 %
From \eqref{:33b} and \eqref{:33c}, we obtain \eqref{:33a}. 
\proofend

\noindent {\em Proof of \tref{l:21}. }
Let $ \mathbb{A}=  
\AAAA _1\ts \cdots \ts \AAAA _m \in \AAAA _n $ 
as in \tref{l:21}. Then, 
because $\AAAA _n \in \Dl  $ for all $ n = 1,\ldots,m $, we deduce 
from \eqref{:10a} and \eqref{:31a} that 
\begin{align}\label{:34a} 
\int_{\mathbb{A} } 
 \rhol ^m(\mathbf{x}) \mmm ^{m} & (d\mathbf{x}) 
= \int_{\mathbb{A} } 
 \rho ^m(\mathbf{x}) \mmm ^{m}(d\mathbf{x}) 
\\ \notag 
= &
\int_{\mathbb{A} } 
 \det \big[
\sum_{\ii ,\jj \in \IIIl } \KFl (\ii ,\jj ) 
\flii (x_p) 
\overline{\fljj (x_q)} \big]_{p,q=1}^m 
\mmm ^{m}(d\mathbf{x}) 
,\end{align}
where $ \mathbf{x}=(x_1,\ldots,x_m)$. 
From a straightforward calculation and \lref{l:31}, we obtain 
\begin{align}\label{:34p}&
\int_{\mathbb{A}  } 
\det \big[
\sum_{\ii ,\jj \in \IIIl } \KFl (\ii ,\jj ) 
\flii (x_p) 
\overline{\fljj (x_q)} \big]_{p,q=1}^m
\mmm ^{m}(d\mathbf{x}) \quad 
\\ 
=\notag & 
\int_{\mathbb{A}  } 
\sum_{\sigma \in \mathfrak{S}_m } 
\mathrm{sgn} (\sigma )
\prod_{p=1}^m 
\Big(
\sum_{\iip , \jjp \in \IIIl } 
 \KFl (\iip ,\jjp ) 
\fliip (x_p) 
 \overline{\fljjp (x_{\sigma (p)})} 
 \Big)
\mmm ^{m}(d\mathbf{x}) 
\\ 
=\notag & 
\sum_{\sigma \in \mathfrak{S}_m } 
\mathrm{sgn} (\sigma )
\int_{\mathbb{A}  } 
\prod_{p=1}^m 
\Big(
\sum_{\iip , \jjp \in \IIIl } 
 \KFl (\iip ,\jjp ) 
\fliip (x_p) 
 \overline{\fljjp (x_{\sigma (p)})} 
 \Big)
\mmm ^{m}(d\mathbf{x}) 
\\ =\notag & 
\sum_{\sigma \in \mathfrak{S}_m } 
\mathrm{sgn} (\sigma )
\limi{R}
\int_{\mathbb{A}  } 
\prod_{p=1}^m 
\Big(
\sum_{\iip , \jjp \in \III (\ell ; R )  } 
 \KFl (\iip ,\jjp ) 
\fliip (x_p) 
\overline{\fljjp (x_{\sigma (p)})}
 \Big)
\mmm ^{m}(d\mathbf{x}) 
\\ =\notag & 
\sum_{\sigma \in \mathfrak{S}_m } 
\mathrm{sgn} (\sigma )
\limi{R}
\int_{\mathbb{A}  } 
\Big(
\sum_{\iii ,\,  \jjj \in \III (\ell ; R ) ^m   } 
\prod_{p=1}^m 
 \KFl (\iip ,\jjp ) 
\fliip (x_p) 
\overline{\fljjp (x_{\sigma (p)})}
 \Big)
\mmm ^{m}(d\mathbf{x}) 
,\end{align}
where $ \III (\ell ; R ) = \{ i \in \IIIl ; \mathrm{rank}(i) \le R \}$ 
and 
$ \mathrm{rank}(i)$ is defined before \eqref{:21j}.
Furthermore, $ \iii = (\ii _1,\ldots,\ii _m) , 
\jjj = (\jj _1,\ldots,\jj _m) \in \IIIl ^m $. 
We note that $ \cup_{i=1}^m \Ai $ is relatively compact. 
Hence the fourth line in \eqref{:34p} follows from 
\lref{l:31} \thetag{2} and the Schwarz inequality. 
Using \lref{l:33}  we obtain 
\begin{align}\label{:34t}&
\int_{\mathbb{A}  } 
\Big(
\sum_{\iii ,\,  \jjj \in \III (\ell ; R ) ^m   } 
\prod_{p=1}^m 
 \KFl (\iip ,\jjp ) 
\fliip (x_p) 
\overline{\fljjp (x_{\sigma (p)})}
 \Big)
\mmm ^{m}(d\mathbf{x}) 
\\ \notag = &
\int_{\mathbb{A}  } 
\Big(
\sum_{\iii ,\,  \jjj \in \III (\ell ; R ) ^m   } 
\prod_{p=1}^m 
 \KFl (\iip ,\jjp ) 
\fliip (x_p) 
\overline{f_{\ell ,\jj _{\sigma^{-1} (p )}} (x_p) }
 \Big)
\mmm ^{m}(d\mathbf{x}) 
\\ \notag 
= &
\int_{\mathbb{A}  } 
\Big(
\sum_{\iii \in \III (\ell ; R ) ^m   } 
\prod_{p=1}^m 
 \KFl (\iip ,\ii _{\sigma (p)} ) 
|\fliip (x_p) |^2
 \Big)
\mmm ^{m}(d\mathbf{x}) 
\\ \notag \to &
\int_{\mathbb{A}  } 
\Big(
\sum_{\iii \in \IIIl ^m   } 
\prod_{p=1}^m 
 \KFl (\iip ,\ii _{\sigma (p)} ) 
|\fliip (x_p) |^2
 \Big)
\mmm ^{m}(d\mathbf{x}) 
\quad \text{ as } R \to \infty 
.\end{align}
The convergence in the last line follows from 
\lref{l:31} \thetag{2} and the Schwarz inequality again. 
Multipling $ \mathrm{sgn} (\sigma )$ and taking summation over 
$ \sigma \in \mathfrak{S}_m $ in the last line, we deduce from 
 \eqref{:21w}--\eqref{:21y} that 
\begin{align}\label{:34r} &
\sum_{\sigma \in \mathfrak{S}_m } \mathrm{sgn} (\sigma )
\int_{\mathbb{A}  } 
\Big(
\sum_{\iii \in \IIIl ^m   } 
\prod_{p=1}^m 
 \KFl (\iip ,\ii _{\sigma (p)} ) 
|\fliip (x_p) |^2
 \Big)
\mmm ^{m}(d\mathbf{x}) 
\\ = \notag & 
\int_{\mathbb{A}  } 
\sum_{\iii  \in \IIIl ^m }
\det [ \KFl (\iip ,\iiq ) 
\big]_{p,q=1}^m 
\Big\{ \prod_{p=1}^m |\fliip (x_p) |^2 \Big\} 
\mmm ^{m}(d\mathbf{x}) 
\\ = \notag &
\int_{\mathbb{A}  } 
\sum_{\iii  \in \IIIl ^m } 
\rhoFl ^m (\iii ) 
\mmFiii (d\mathbf{x}) 
\\ = \notag &
\sum_{\iii  \in  \in \III_{\ell }(\mathbb{A}) }
\rhoFl ^m (\iii ) 
.\end{align}
Combining \eqref{:34a}--\eqref{:34r} 
we deduce \eqref{:21b}, which completes the proof.  
\qed

\section{Proof of \tref{l:22}--\tref{l:25}}\label{s:4} 

In this section, we prove \tref{l:22}--\tref{l:25}. 

\subsection{Proof of \tref{l:22}} 
Let $ \varrho ^m $ be the $ m $-point correlation function of $\8 $. 
Then it suffices for \eqref{:22a} to prove 
\begin{align}\label{:41z}& 
\rhol ^m (\mathbf{x}) = \varrho ^m (\mathbf{x})  
.\end{align}

From \eqref{:20e} and $ \mathcal{F}\kl = \sigma [ \Ali ; i\in \Il  ]$, 
we see that $ \rhol ^m $ and $ \varrho ^m $ are $ \Fl ^m $-measurable. 
Let $  m = m_1+\cdots + m_k $. Let 
$\mathbb{A} = \AAAA _1^{m_1}\ts \cdots \ts \AAAA _k^{m_k} \in\Dl ^m $ 
such that 
$ \AAAA _p \cap \AAAA _q = \emptyset $ if $ p\not=q $. 
Let 
$ \mathbf{i} = (i_n)_{n=1}^m = 
(\mathbf{i}_1,\ldots,\mathbf{i}_k)\in \IIIl ^m  $ 
such that $ \mathbf{i}_n \in \III (\ell )^{m_n} $. 
From \tref{l:21}, we see that 
\begin{align}\label{:41a}
\int_{\mathbb{A}} \rhol ^m(\mathbf{x}) \mmm ^{m}(d\mathbf{x}) = 
\sum_{\iii \in \III_{\ell }(\mathbb{A}) } \rhoFl ^m (\iii ) 
.\end{align}
By the definition of correlation functions, 
\eqref{:22w}, and \eqref{:22x}, we see that 
\begin{align}\label{:41d}
\sum_{\iii \in \III_{\ell }(\mathbb{A}) } 
\0  
=&
\int _{\IIl }\prod_{n=1}^k
\frac{\mathsf{i}  (\III_{\ell }(\AAAA _n)) ! }
{(\mathsf{i}  (\III_{\ell }(\AAAA _n)) -m_n)!} 
\nul (d\mathsf{i})
\\ \notag 
=& 
\int_{\Ssf }
\prod_{n=1}^k
\frac{\mathsf{s} (\AAAA _n ) ! }{(\mathsf{s} (\AAAA _n ) -m_n)!}  
\8 (d\mathsf{s})
\\ \notag 
=&
\int_{\mathbb{A}} 
\varrho ^m (\mathbf{x}) \mmm ^m (d\mathbf{x})
.\end{align}
Combining \eqref{:41a} and \eqref{:41d}, we deduce that 
\begin{align}\label{:41e}&
\int_{\mathbb{A}} 
\rhol ^m (\mathbf{x}) \mmm ^m (d\mathbf{x}) =
\sum_{\iii \in \III_{\ell }(\mathbb{A}) } \0  
= 
\int_{\mathbb{A}} 
\varrho ^m (\mathbf{x}) \mmm ^m (d\mathbf{x})
.\end{align}
From \eqref{:41e}, we obtain \eqref{:41z}. This completes the proof 
of \tref{l:22}. 
\qed

\subsection {Proof of \tref{l:23}} 
\tref{l:23} follows from \tref{l:22} immediately. 
\qed

\subsection{Proof of \tref{l:24}}

It is known that determinantal point processes on discrete spaces are tail trivial 
\cite{RL.03,s-t.jfa}. 
Hence $ \nuFl  $ is tail trivial by \lref{l:32}. 

Let $ \ulabl $ be as in \tref{l:22}. 
Let $ \mathsf{A} \in \ulabl ^{-1} (\Gl ) $. 
Then there exists a $ \mathsf{B} \in \mathcal{B}(\IIl ) $ such that 
$ \mathsf{A} = \iotal ^{-1} (\mathsf{B})$. 
From this we deduce that, for each $ \mathsf{A} \in \Tail (\4 ) \cap \ulabl ^{-1} (\Gl ) $, 
there exists 
a $ \mathsf{B} \in \Tail (\IIl ) $ such that 
$ \mathsf{A} = \iotal ^{-1} (\mathsf{B})$. 
Hence from \eqref{:22w} we deduce 
\begin{align}\label{:42q}&
\7 (\mathsf{A}) = \nuFl (\mathsf{B})
.\end{align}
From \eqref{:42q} and tail triviality of $ \nuFl  $ we deduce that 
\begin{align}\label{:42r}&
\7 (\mathsf{A}) \in \{ 0,1 \} 
\end{align}
for each $ \mathsf{A} \in \Tail (\4 ) \cap \ulabl ^{-1} (\Gl ) $. 
We easily see that 
$ \ulabl ^{-1} (\Gl ) \subset \sigma [\iotal ] $.  
Hence 
\begin{align}\label{:42p}&
 \Tail (\4 ) \cap \ulabl ^{-1} (\Gl ) \subset 
\Tail (\4 ) \cap \sigma [\iotal ] 
.\end{align}

Combining \eqref{:42r} and \eqref{:42p} completes the proof of \tref{l:24}.  
\qed

\subsection{Proof of \tref{l:25}} 

Let $ \mathsf{B} \in \TailS \cap \Gl $. Then we deduce that
\begin{align}\notag &
 \ulabl ^{-1} (\mathsf{B}) \in \Tail (\4 ) \cap \ulabl ^{-1} (\Gl ) 
.\end{align}
Hence from \tref{l:22} and \tref{l:24}, we deduce that 
\begin{align}\notag &
\mu (\mathsf{B}) = \mul (\mathsf{B}) 
= \7 (\ulabl ^{-1} (\mathsf{B}) ) \in \{ 0,1 \} 
.\end{align}
This completes the proof. 
\qed 

\section{Proof of \tref{l:11}}\label{s:6}

In this section, we complete the proof of \tref{l:11}.  
%

\begin{lemma} \label{l:62}
Let $ X $ be a $ \Tail (\Ssf )$-measurable and integrable random variable. 
Then $ E^{\mu }[X|\Gl ]$ is $ \TailS \cap \Gl $-measurable. 
\end{lemma}

\proofbegin  
Recall that $ \Dl = \{ \Ali \}_{i\in\Il } $. 
Let $ \pi _{T_r}$ be the projection with $ T_r $ such that 
\begin{align}\label{:62a}&
T_r = \bigcup_{\Ali \cap \Sr \not= \emptyset ; \atop i\in \Il } 
\Ali 
.\end{align}
Then $ X \in L^1(\Ssf ,\mu ) $ is $ \sigma [{\pi _{T_r^c}}] $-measurable 
because $ X \in L^1(\Ssf ,\mu ) $ is $ \Tail (\Ssf )$-measurable and 
each $ \Ali $ is relatively compact. Hence for each $ r \in \N $
\begin{align}\label{:62b}&
X (\mathsf{s}) 
 = X \circ \pi _{T_r^c} (\mathsf{s}) 
.\end{align}
From this we deduce that 
\begin{align}\label{:62c}&
E^{\mu } [X | \Gl ] = E^{\mu } [ X \circ \pi _{T_r^c} | \Gl ] 
.\end{align}
By construction $ \Sr \subset T_r $. 
Then from this and \eqref{:62c} we see that 
$ E^{\mu } [X | \Gl ] $ is 
$ \sigma [\pi _{S_r^c} ]$-measurable for each $ r \in \N $. 
Hence  
$ E^{\mu }[X|\Gl ]$ is $ \TailS $-measurable 
because 
$ \cap_{r\in\N } \sigma [\pi _{S_r^c} ]=\Tail (\Ssf )$. 
By construction  $ E^{\mu } [X | \Gl ] $ is 
$ \cap_{r\in\N } \sigma [\pi _{S_r^c} ]$-measurable. 
Combining these completes the proof of \lref{l:62} 
\proofend

\begin{lemma} \label{l:63} For all $ \mathsf{A} \in \TailS $ 
\begin{align}\label{:63a}&
\mu(\mathsf{A}) = \mulsA  
\quad \text{ for $ \mu $-a.s.\! $ \mathsf{s}$}
.\end{align}
\end{lemma}
\proofbegin 
From the definition of the conditional probability, we see that 
\begin{align}\label{:63b}&
\mu (\mathsf{A}) = \int_{\Ssf } \mulsA  \mu (d\mathsf{s})
.\end{align}
From \lref{l:62},  
we deduce that 
$ \mulsA  = E^{\mu}[1_{\mathsf{A}}|\Gl ](\mathsf{s})$ is 
$ \TailS \cap \Gl $-measurable. 
Hence from \tref{l:25} we obtain that $ \mulsA  $ is constant 
$ \mu $-a.s.\! $ \mathsf{s}$. 
This combined with \eqref{:63b} yields \eqref{:63a}.  
\proofend

\begin{lemma} \label{l:61}
For each $ \mathsf{A} \in \mathcal{B}(\Ssf )$ 
\begin{align}\label{:61z}& \quad \quad 
\limi{\ell } \mulsA  = 
1_{\mathsf{A}} (\mathsf{s}) 
\quad \text{ for $ \mu $-a.s.\!\! $ \mathsf{s}$}
.\end{align}
\end{lemma}
\proofbegin 
From \eqref{:20f}, we apply the martingale convergence theorem 
to obtain the convergence such that, for all $ \mathsf{A} \in \mathcal{B}(\Ssf ) $, 
\begin{align}\label{:61b}& 
\limi{\ell } \mulsA  = 
\limi{ \ell } 
E^{\muk }[1_{\mathsf{A}}| \Gl ] (\mathsf{s}) = 
E^{\muk } [1_{\mathsf{A}}|\mathcal{B}(\Ssf ) ]  (\mathsf{s}) = 
1_{\mathsf{A}} (\mathsf{s}) 
\end{align}
 for $\muk $-a.s.\! $ \mathsf{s}$.. We have thus proved \eqref{:61z}.  
\proofend

\noindent 
{\em Proof of \tref{l:11}.} 
From \lref{l:63} and \lref{l:61} we deduce that 
\begin{align}\label{:64b}&
\mu (\mathsf{A}) =  \mulsA \to _{\ell \to \infty }1_{\mathsf{A}} (\mathsf{s}) 
\quad \text{ $ \mu $-a.s.\! $ \mathsf{s}$} 
.\end{align}
Hence we obtain $ \mu (\mathsf{A}) \in \{ 0,1 \} $. 
\qed

\section{Examples related to random matrices} \label{s:7}

In this section, we give typical examples of determinantal point processes 
related to random matrix theory \cite{Meh04,forrester}. All examples below are tail trivial because of \tref{l:11}.  

All the kernels $ \K (x,y)$ below are continuous. 
In Examples \ref{d:7sin}--\ref{d:7be}, 
we define the kernels only off diagonal. 
On diagonal, they are defined by continuity. 

\begin{example}[sine point process] 
\label{d:7sin} 
Let $ \sS = \mathbb{R}$ and $ \mmm (dx) = dx $. Let 
\begin{align} \notag 
&
\K _{\mathrm{sin}}(x,y) = \frac{\sin (x-y)}{\pi (x-y) } \quad (x\not=y)
\end{align}
be the sine kernel. 
The associated determinantal point process  
$ \mu _{\mathrm{sin}}$ is called the sine$_{2}$ point process.

\end{example}

\begin{example} [Airy point process]\label{d:76} 
Let $ \sS = \mathbb{R}$ and $ \mmm (dx) = dx $. Let 
\begin{align} \notag 
&
 \mathsf{K}_{\mathrm{Ai}}(x,y) = 
 \frac{\mathrm{Ai}(x)\mathrm{Ai}'(y) - \mathrm{Ai}'(x)\mathrm{Ai}(y) }{x-y} 
 \quad (x\not=y)
\end{align}
be the Airy kernel. Here 
$ \mathrm{Ai}$ is the Airy function, and $ \mathrm{Ai}'$ is its derivative. 
The associated determinantal point process $ \mu _{\mathrm{Ai}}$ is called 
the Airy point process \cite{Meh04,forrester}.  
\end{example}

\begin{example}[Bessel point process]\label{d:7be} 
Let $ \sS = [0,\infty) $ and $ \mmm (dx) = dx $. 
Let $ 1 \le \alpha < \infty $. 
Let $ \mathsf{K}_{\mathrm{Be},\alpha } $ be the Bessel kernel such that  
\begin{align}\notag 
&
 \mathsf{K}_{\mathrm{Be},\alpha }(x,y) = 
 \frac{J_{\alpha } (\sqrt{x}) \sqrt{y} J_{\alpha }' (\sqrt{y}) - 
 \sqrt{x} J_{\alpha }' (\sqrt{x}) \sqrt{y} J_{\alpha }(\sqrt{y}) }{2(x-y)}
\quad (x\not=y)
.\end{align}
Let $ \mu _{\mathrm{Be},\alpha }$ 
be the associated determinantal point process.  
$ \mu _{\mathrm{Be},\alpha }$ is called the Bessel $_{2,\alpha }$ point process. 
\end{example}

\smallskip 

\begin{example}[Ginibre point process] \label{d:74} 
Let $ \sS = \mathbb{R}^2$ and $ \mmm (dx) = (1/\pi) e^{-|x|^2}dx $. 
Let $ \map{\kg }{\mathbb{R}^{2}\ts \mathbb{R}^{2}}{\mathbb{C}}$ 
be the exponential kernel such that 
\begin{align} \notag 
&
\kg (x,y) = 
e^{x \bar{y}}
.\end{align}
Here we identify $ \mathbb{R}^{2} $ as $ \mathbb{C}$ by the obvious correspondence 
$ \mathbb{R}^{2} \ni x=(x_1,x_2)\mapsto x_1 + \sqrt{-1}  x_2 \in \mathbb{C}$, and 
$ \bar{y}=y_1-\sqrt{-1}  y_2 $ is the complex conjugate in this identification. %
The associated determinantal point process $ \mug $ is called the Ginibre point process. 
\end{example}


\bs
\bs

\noindent 
{\bf Acknowledgement: }\\
H.O. is supported in part by a Grant-in-Aid for Scientific Research
(KIBAN-A, No.\! 24244010)
from the Japan Society for the Promotion of Science.
\end{document}